\documentclass[10pt]{article}
\usepackage{amsmath, verbatim}
\usepackage{graphicx,color}
\newcommand{\qed}{{\unskip\nobreak\hfil\penalty50\hskip2em\vadjust{}
            \nobreak\hfil$\Box$\parfillskip=0pt\finalhyphendemerits=0\par}}

\newtheorem{thm}{Theorem}[section] 
\newtheorem{lemma}{Lemma}[section] 

\newtheorem{definition}{Definition}[section]
\newcommand{\bed}{\begin{definition}}
\newcommand{\eed}{\end{definition}}

\newcommand{\eps}{\epsilon}

\newcommand{\bitem}{\begin{itemize}}
\newcommand{\eitem}{\end{itemize}}

\newcommand{\goto}{\rightarrow}

\newcommand{\beqn}{\begin{equation}}
\newcommand{\eeqn}{\end{equation}}
\newcommand{\balign}{\begin{align}}
\newcommand{\ealign}{\end{align}}

\newcommand{\so}{\sigma_0}
\newcommand{\uo}{u_0}

\usepackage{amssymb}
\begin{document}
\nocite{*}
\title{A Generalized Fourier Approach to Estimating the Null Parameters and Proportion of nonnull Effects  in Large-Scale Multiple Testing}
\author{Jiashun Jin$^1$ and Jie Peng$^2$ and Pei Wang$^3$ \\
$^1$Statistics Department,  Carnegie Mellon  University\\
$^2$Statistics Department, University of California at Davis \\
$^3$Division of Public Health Sciences, Fred Hutchinson Cancer
Research Center} \maketitle
\begin{abstract}
In a recent paper  \cite{Efron},  Efron pointed out that an
important issue in large-scale multiple hypothesis testing is that
the null distribution may be unknown and need to be estimated.
Consider a Gaussian mixture model, where the null distribution is
known to be normal but both null parameters---the mean   and the
variance---are unknown.  We address the problem with a method based
on Fourier transformation.   The Fourier approach was first studied
by Jin and Cai \cite{JC}, which focuses on the scenario where any
non-null effect has either the same or a larger variance than that
of the null effects. In this paper, we review the main ideas in
\cite{JC}, and propose a generalized Fourier approach to tackle the
problem under another scenario: any non-null effect has a larger
mean than that of the null effects, but no constraint is imposed on
the variance. This approach and that in \cite{JC} complement with
each other: each approach is successful in a wide class of
situations where the other fails. Also,   we extend the Fourier
approach to estimate the proportion of non-null effects. The
proposed procedures perform well both in theory and on simulated
data.
\end{abstract}

\begin{quote} \small
{\bf Keywords}:   empirical null,  Fourier transformation,  generalized Fourier transformation,  proportion of non-null effects,  sample size calculation,
\end{quote}
\begin{quote} \small
{\bf AMS 1991 subject classifications}:   Primary  62G10, 62G05;  secondary 62H15, 62H20.
\end{quote}
{\bf Acknowledgments}: The authors are in part supported by grants
from National Science Foundation, DMS-0908613 (Jin), DMS-0806128 (Peng) and a grant from
National Institute of Health, R01GM082802 (Peng and Wang).

\section{Introduction}
Large-scale multiple testing is a recent  area of  active research  in statistics, where one tests  thousands or  even millions of null hypotheses {\it simultaneously}:
\[
H_j,  \qquad   j = 1, \ldots, n.
\]
Associated with each null hypothesis is a test statistics $X_j$,
which, depending on the situation, can be a summary statistic, a
$p$-value,  a regression coefficient, or a transform coefficient,
etc..     We say that  $X_j$ contains a {\it null effect}  if $H_j$
is true, and contains a {\it non-null effect} if  otherwise.

A convenient model is the  Bayesian hierarchical model \cite{Efron,
Wasserman} which we now describe. Fix $0 < \eps  < 1$. For each $1
\leq j \leq n$,   we flip a coin with probability $\eps$  of landing
tail.   If the coin lands head, we draw $X_j$ from a common density
function $f_0(x)$  which we call the {\it null density}. If the coin
lands tail,      we draw $X_j$ from an individual density function
$\xi_j(x)$,  where $\xi_j$ itself is randomly generated according to
a fixed probability measure  $\Xi$. In effect, $X_j$ can be viewed
as samples from the density  $f_1(x) \equiv \int \xi(x) d \Xi(\xi)$,
which we call the {\it alternative density}; see \cite{Wasserman,
Jin}.
 Marginally,    $X_j$ can be deemed as samples from the following two-component mixing density:
\begin{equation}  \label{Eqmodel1}
X_j \;\;  \stackrel{iid}{\sim}   \;\;    (1 - \eps)   f_0(x)  + \eps   f_1(x)  \;\;   \equiv \;\;   f(x).
\end{equation}

The parameter $\eps$  is closely related to the proportion of
non-null effects (i.e., the fraction of null hypotheses that are
untrue).  In fact, under the Gaussian mixture model, the number of
untrue hypothesis is distributed as Binomial with parameters n and
$\eps$. So when n is large, the difference between $\eps$ and the
actual fraction $\leq O_p(\sqrt{\eps/n})$ and is usually
negligible. For this reason, we call $\eps$ the proportion of the
non-null effects in this paper.

The null density is the starting point for any testing procedures.
In many scenarios, the null density is assumed as known.  However,
somewhat surprisingly, this assumption may be incorrect in some
multiple testing situations as pointed out by Efron \cite{Efron}.
Efron illustrated his point with a breast cancer microarray data,
which is based on $15$ patients with $7$ having BRCA1 mutation and
$8$ having BRCA2 mutation. For each patient, the same set of $3226$
genes were measured and it is of interest to find which genes are
differentially expressed.  For each gene, a studentized-$t$ score
was calculated and then transformed to a $z$-score (see
\cite{Efron}) for the details). Efron argued that, although the
theoretical  null should be the standard normal $N(0,1)$,   another
null density, $N(0.02, 2.50)$ seems to be more appropriate.  Efron
called the later the  {\it empirical null} and demonstrated
convincingly  that it is better to use the empirical null  instead
of the theoretical null in many situations.

There are many possible reasons
why the empirical null may be different from the theoretical  null.
 Take the breast  cancer microarray data for example, the studentized-$t$ statistics may not be truly $t$-distributed due
to failed  distributional assumptions. There may be  covariates
(such as age  of the patients) that has not been observed  in the
data.   The correlation across different genes (also that across
different arrays)  has   been neglected.   All these factors may
drive the empirical null far   from the theoretical  null.

Unfortunately, unlike the theoretical  null, the empirical null is
usually unknown.  Thus how to estimate the empirical null is
a problem of major interest.

\subsection{Identifiability issue and constrained Gaussian mixture models}
Note that in Model (\ref{Eqmodel1}), some may call $f_0$ the null
density, and some may call $f_1$ the null density. To resolve this
issue, we fix a constant $\eps_0 \in (0,1/2)$ and  assume
\[
0 < \eps  \leq \eps_0,
\]
so that the null density is tied to the majority of the hypotheses.

We adopt the Gaussian model as suggested in  Efron \cite{Efron}. In
detail,  let  $\phi(\cdot)$ be the density of  $N(0,1)$.
We assume  that the null density $f_0$ is Gaussian   with an unknown
mean $\uo$ and an unknown variance $\so^2$:
\[
f_0(x) =  \frac{1}{\so}  \phi(\frac{x - \uo}{\so}).
\]
At the same time, we assume that  the alternative density $f_1$ is a Gaussian mixture (both a location mixture and a scale mixture)  with a bivariate mixing distribution  $H(u, \sigma)$:
\[
f_1(x)  =   \int \frac{1}{\sigma} \phi(\frac{x - u}{\sigma}) d H(u, \sigma).
\]
The marginal density of $X_j$ is then
\begin{equation} \label{Definef}
f(x)  = f(x;  \uo, \so,  \eps,  H) =   (1 - \eps)  \frac{1}{\so} \phi(\frac{x - \uo}{\so})  +  \eps   \int \frac{1}{\sigma} \phi(\frac{x - u}{\sigma}) d H(u, \sigma).
\end{equation}
With the  Gaussian model, the problem of  estimating the null density  reduces  into the problem of estimating the null parameters $(\uo, \so^2)$.

However, the null density in the above Gaussian model is not always
identifiable. This is because, without constraint on $H(\cdot, \cdot)$,  $f_1$
can be very close or even identical to $f_0$. Fortunately, there are many natural
constraints that we can put on $H(\cdot, \cdot)$ to resolve this
problem. Below are some examples.

\bed Fix $\eps_0 \in (0,1/2)$,  $\uo$, and $\so > 0$.  We say that
$f(x) = f(x; \uo, \so, \eps,  H)$ is a Gaussian mixture density
constrained with Elevated Variances with respect to parameters
$(\uo, \so,\eps_0)$    if it has
the form as in (\ref{Definef}),  and that   the proportion $\eps$
and the mixing distribution $H$ satisfy
\begin{equation} \label{GEV}
0 < \eps  \leq \eps_0, \qquad  P_{H}  (\sigma \geq  \so) = 1,  \qquad P_{H} \bigl((u, \sigma)   \neq      (\uo, \so) \bigr)  =  1.
\end{equation}
We refer to the Gaussian model (\ref{Definef}) with constraints in
(\ref{GEV}) as $GEV(\uo,
\so, \eps_0)$. \eed \vspace{-5pt} For short, we write $GEV(\uo,
\so, \eps_0)$ as $GEV$ whenever there is no confusion. In the
definition above,  $u$ and $\sigma$ denote the location and scale
parameters from the mixing distribution, and $P_{H}$ denotes the
probability under the mixing distribution $H(\cdot, \cdot)$. This
models a situation where the variance associated with an individual
non-null effect is no less than that of a null effect. The following
lemma shows that, given (\ref{GEV}), the triplets $(\uo, \so, \eps)$
are uniquely determined by $f(x)$ and the identifiability issue is
therefore resolved. \vspace{-4pt}
\begin{lemma} \label{lemma:GEV}
Given a density $f(x) = f(x;  \uo, \so, \eps,  H)$ satisfying  (\ref{Definef}) and  (\ref{GEV}),  the parameters $\uo$, $\so$, and $\eps$ are uniquely determined by $f(x)$.
\end{lemma}\vspace{-4pt}
This lemma is proved in Section \ref{sec:appen}. Note that if  we replace the
constraint  $P_{H} (\sigma \geq \so) = 1$  by  $P(\sigma
\leq \so) = 1$, then the identifiability issue persists (the
construction of counter examples  is elementary and we skip
it).\vspace{7pt}

 Alternatively, we define GEM   as follows. \vspace{-5pt}
\bed Fix $\eps_0 \in (0,1/2)$,  $\uo$, and $\so > 0$.  We say that
$f(x) = f(x; \eps, \uo, \so, H)$ is a Gaussian mixture  density
constrained with Elevated Means with respect to parameters $(\uo,
\so, \eps_0)$    if it has
the form as in (\ref{Definef}),  and  that  the proportion $\eps$
and the mixing distribution $H$ satisfy
\begin{equation} \label{GEM}
0 < \eps  \leq \eps_0, \qquad  P_{H}  (u  >    \uo) = 1.
\end{equation}
We refer to the Gaussian model (\ref{Definef}) with constraints in
(\ref{GEM}) as  $GEM(\uo, \so, \eps_0)$. \eed \vspace{-3pt} GEM models a situation where
the mean associated with an individual non-null effect is larger than that of a null effect. The following lemma is proved in Section \ref{sec:appen}.\vspace{-4pt}
\begin{lemma} \label{lemma:GEM}
Given a density $f(x) = f(x;  \uo, \so, \eps,  H)$ satisfying  (\ref{Definef}) and  (\ref{GEM}),   the parameters $\uo$, $\so$, and $\eps$  are uniquely determined by $f(x)$.
\end{lemma}\vspace{-4pt}
For the case where we replace the  constraint $P_{H}(u >   \uo) = 1$
in (\ref{GEM}) with $P_{H}(u  <    \uo)  = 1$,   the discussion is
similar.  Also,  we can relax the constraint to $P_{H} (u \geq \uo)
= 1$. But by doing so we need some conditions on  $\sigma$.  For reasons of space, we skip the discussion along these
two lines.

GEV and GEM are the two main models we study in this paper. Despite
the additional constraints, both models are broad enough to
accommodate  many interesting cases that arise in real applications.
In sections below, we discuss possible approaches to consistently
estimating the null parameters in GEV and GEM.

\subsection{A Fourier approach to estimating the null parameters in GEV}
Conventionally, one estimates the null parameters with either
empirical moments or extreme observations.  However,  in these
quantities,  the information containing the null parameters is
highly distorted by the non-null effects. A non-orthodox approach is
therefore necessary. In a recent work \cite{JC},  Jin and Cai
proposed a Fourier approach to estimating the null parameters in
GEV. We now briefly explain the idea.

When it comes to a density function, one usually  pictures it as a smooth curve that spreads over the real line.  Joseph Fourier taught us a different view point:  a normal density $N(u, \sigma^2)$ is not only a bell shaped
curve centered at $u$, but also a wave oscillate at the frequency $u$. In fact, the Fourier transform of the density   $N(u, \sigma^2)$ can be decomposed into two components:  the amplitude function  determined by $\sigma^2$, and the phase function  determined by $u$:
\begin{equation} \label{Phase1}
e^{-\sigma^2 t^2/2}  \cdot   e^{it u}   \equiv  \mbox{Amplitude} \cdot   \mbox{Phase function}, \qquad i  = \sqrt{-1}.
\end{equation}
Consequently,  we can view a Gaussian mixture as a superposition of waves with different frequencies and different amplitudes.

We now invoke GEV.  The above investigation gives rise to an
interesting approach to estimating the null parameters.  Denote the
{\it empirical characteristic function}  by
\[
\psi_n(t) = \psi_n(t; X_1, \ldots, X_n) =  \frac{1}{n} \sum_{j = 1}^n e^{i t X_j}.
\]
For an appropriately large frequency $t$, the stochastic fluctuation
is negligible and $\psi_n$  reduces to its non-stochastic
counterpart---the {\it underlying characteristic function}
$\psi(t)=E[\psi_n(t)]$. By direct calculations,
\[
\psi(t) =  \psi(t; \uo, \so, \eps, H)  \equiv  \psi_0(t) [1 + s(t)],
\]
where
\[
\psi_0(t) = \psi_0(t; \uo, \so, \eps) = (1 - \eps) e^{i \uo t - \so^2 t^2/2},
\]
and
\begin{equation}
\label{Definess}
s(t) =  s(t; \uo, \so, \eps, H) =    \frac{\eps}{1 - \eps}  \int   e^{i  (u - \uo)  t -  (\sigma^2 - \so^2)  t^2/2} d H(u, \sigma).
\end{equation}

Now,  with  GEV and a little bit extra condition,    $s(t) \approx
0$.  For example, if we assume that    $P_{H}(\sigma > \so) = 1$,
then at a high frequency $t$,
\begin{equation} \label{Defines}
|s(t)| \leq \frac{\eps}{1 - \eps}  \int   e^{-  (\sigma^2 - \so^2)  t^2/2} d H(u, \sigma)  \approx 0.
\end{equation}
This says that in GEV,  as the frequency $t$ tends to $\infty$, the waves corresponding to the alternative density damps faster than that associated with
the null density.  Therefore,  the information containing  the null parameters is
asymptotically preserved in high frequency Fourier transform, where the distortion of non-null effects is negligible. In other words,   for an appropriately large frequency $t$,
\[
\psi_n(t) \approx \psi(t) \approx \psi_0(t).
\]
Now, since $\psi_0(t)$ has a very simple form,   we can solve  $(\uo, \so)$ (and also $\eps$)  from it.

The elaboration of the idea gives rise to the estimators in \cite{JC}, which are proved to be uniformly consistent to the  null parameters across a wide range  of mixing
distributions  $H(\cdot, \cdot)$. It was also shown in \cite{JC} that these  estimators  attain the optimal rate of convergence. See the details therein.    These works reveal that, somewhat surprisingly,   the right place to estimate the null parameters is in the frequency domain, rather than in  the spatial domain as one may have expected.

\subsection{A generalized-Fourier approach to estimating the null parameters in GEM}
Despite its encouraging performance  in GEV,  the above approach does not yield a satisfactory estimation in GEM.  To see the point, we  note that  the key for the success of  the above approach is (\ref{Defines}),   which critically depends on the assumption of $P_{H} (\sigma > \so) = 1$.  Note that such an assumption does not hold in GEM.  As a result,    the above approach  ceases to perform well.

Fortunately, there is an easy fix.     The key is to  replace the
Fourier transformation by the generalized Fourier transformation (to
be introduced below),  so that  in the frequency domain, the roles of
the mean and the variance are  sort of ``swapped".   In
detail, let
\[
\omega  =   - (1+i)/\sqrt{2},     \qquad (\mbox{ note $\omega^2 = i$}).
\]
For any density function $h(x)$, the generalized Fourier transformation is
\[
\int h(x) \exp(\omega x) dx,
\]
provided that the function $h(x) \mathrm{exp}(\omega x)$ is absolutely  integrable.  In particular, the generalized Fourier transform of the Gaussian density $N(u, \sigma^2)$ is
\begin{equation} \label{Phase2}
\exp \bigl( - \frac{ut }{\sqrt{2}} \bigr)   \cdot  \exp \bigl( i  [ -  \frac{ut }{\sqrt{2}} + \frac{\sigma^2t^2}{2}] \bigr) \equiv \mbox{Amplitude function} \cdot \mbox{Phase function}.
\end{equation}
Now,  the amplitude is uniquely determined by the mean  (compare with (\ref{Phase1})).

The remaining part of the idea is similar to that in the preceding section.
Denote the {\it generalized-empirical characteristic function}  by
\[
\varphi_n(t) =  \varphi_n(t; X_1, \ldots, X_n) = \frac{1}{n} \sum_{j = 1}^n \exp( \omega t X_j).
\]
For large $n$ and an appropriately chosen $t$, one expects that  the stochastic fluctuation is negligible, and that $\varphi_n(t)$ reduces approximately to the {\it generalized characteristic function},
\[
\varphi(t)  = \varphi(t; \uo, \so,  \eps, H)   \equiv E[\varphi_n(t)].
\]
Direct calculations show that
\[
\varphi(t) =       \varphi_0(t)  [1 + r(t)],
\]
where
\[
\varphi_0(t) = \varphi_0(t; \uo, \so, \eps)  = (1 - \eps)    \exp ( \omega  u_0 t  +   i  \sigma_0^2 t^2 /2),
\]
and
\begin{equation} \label{Definer}
r(t) = r(t; \uo, \so, \eps, H)   =  \frac{\eps}{1-\eps}\int  \exp( \omega(u - \uo) t  + i (\sigma^2 - \so^2)t^2/2) d H(u, \sigma).
\end{equation}
Recalling that  $\omega =  - (1 + i)/\sqrt{2}$, it is seen that
\[
|r(t)| \leq  \frac{\eps}{1-\eps} \int   \exp(- (u - \uo)t/\sqrt{2}) d H(u, \sigma).
\]

We now invoke GEM.  Similarly,  since that  $P_{H} (u > \uo) = 1$,   $r(t) \approx 0$ for large $t$.
We  expect that
\[
\varphi_n(t)  \approx \varphi(t) \approx \varphi_0(t).
\]
Again,   $\varphi_0(t)$ has a very simple form and we can solve $(\uo, \so)$ (and also $\eps$) from it.   In fact, introduce two functionals
$\uo(\cdot; t)$ and $\so^2(\cdot;t)$ by
\begin{equation} \label{Functional}
\uo(g;t) =-\frac{\sqrt{2}}{|g(t)|}\frac{d}{d t}|g(t)|, \qquad
\sigma^{2}_{0}(g;t) =\frac{\sqrt{2} \mathrm{Re}(\omega \bar{g} g') }{t |g(t)|^2},
\end{equation}
where $g$ is any complex-valued differentiable function,  and $|z|$, $Re(z)$ and $\bar{z}$  denote the module,  the  real part,   and the complex conjugate of  a complex number $z$, correspondingly.
The following lemma says that plugging $g = \varphi_0$ into two functionals gives the desired parameters $\uo$ and $\so^2$, respectively.
\begin{lemma} \label{lemma:varphi0}
For all $t \neq 0$, $\uo(\varphi_0; t)  = \uo$ and $\so^2(\varphi; t) = \so^2$.
\end{lemma}
Lemma \ref{lemma:varphi0}  can be proved using elementary algebra, so we skip it.    Taking $g = \varphi_n$ in  (\ref{Functional}),   we expect to have
\[
\uo(\varphi_n, t) \approx \uo(\varphi, t) \approx \uo, \qquad \so^2(\varphi_n, t) \approx \so^2(\varphi, t) \approx \so^2.
\]

In this paper, we shall carefully study the bias  and variance  of $\uo(\varphi_n; t)$ and $\so^2(\varphi_n;t)$,   and investigate which choices of $t$ give  a   good tradeoff between the bias and the variance.   We find out that as $n$ tends to $\infty$, if we set $t$ in an appropriate range,  then both estimators
are consistent with their estimands,  uniformly so across a wide class of situations.

\subsection{Estimating the proportion of non-null effects}
Seemingly, the approach can be readily generalized to estimate the proportion of  non-null effects $\eps$.
 How to estimate the proportion  has been the topic of
many recent works in
 the area of large-scale multiple hypothesis testing.  See for example \cite{CJ, Efronetal, Wasserman, Jin, JC,  Langaas,  Rice,  Storey1,  Swanepoel}.  There are two
reasons for the enthusiasm.  In some applications, the proportion is
the quantity that is of direct interest \cite{Rice}; while more
often, knowing the proportion helps to improve many multiple testing
procedures, such as the  FDR procedure by Benjamini and Hochberg's
\cite{Benjamini}, the local FDR procedure by Efron et al.
\cite{Efronetal} and  the optimal discovery function by Storey
\cite{Storey2}. See \cite{Jin} for more discussions.

In Section \ref{sec:prop},  we
extend the generalized Fourier
approach to estimating the proportion in GEM.   We discuss  two
different cases: (1) the null parameters are known; and (2) the null
parameters are unknown.  In both cases, we find that the estimators
are uniformly consistent with the proportion across a wide class of
situations.

We remark that the success of  the Fourier approach for estimating
the null parameters and the proportion is not  coincidental.  It
roots from the key fact  that the null density can be isolated from
the alternative density in the high frequency Fourier coefficients.
Naturally,  we shall  continue to find the Fourier approach to be
successful in estimating many other quantities.

The remaining part of the paper is organized as follows.    Section
\ref{sec:main} studies the problem of estimating the null parameters
in GEM. We  show that by choosing an appropriate $t$,  the
estimators $\uo(\varphi_n; t)$ and $\so^2(\varphi_n; t)$ are
consistent to the true parameters, uniformly across a wide class of
situations. Section \ref{sec:prop} studies the problem of estimating
the proportion. While the studies  in Sections
\ref{sec:main}--\ref{sec:prop} are  asymptotic, we  carry out a few
simulation studies in Section \ref{sec:simul}, and investigate the
performance of the proposed estimators for moderately large $n$.
Section \ref{sec:appen} contains the proofs for the theorems and
lemmas, in the order  they  appear.

\section{Main results} \label{sec:main}
In this section,  we limit our attention to GEM and  study the estimation errors of $\uo(\varphi_n; t)$ and $\so^2(\varphi_n; t)$.   Since the discussions are similar, we focus on that of $\uo(\varphi_n; t)$.
For the asymptotic analysis,   we adopt a framework where  both  $\eps$ and $H$  may depend on  $n$ as $n$ ranges from $1$ to $\infty$ (denoted by $\eps_n$ and $H_n$).  This covers a much  broader situations than that  when
$(\eps, H)$ are fixed as  $n$ ranges from $1$ to $\infty$.

\subsection{Asymptotic framework}
Recall that the test statistics $X_j$ are iid samples from
\begin{equation} \label{f}
f(x) = f(x; \uo, \so, \eps_n, H_n, n) = (1 - \eps_n) \frac{1}{\so}  \phi(\frac{x - \uo}{\so}) + \eps_n \int \frac{1}{\sigma} \phi(\frac{x  -  u}{\sigma}) d H_n(u, \sigma).
\end{equation}
As before,   fix $\eps_0 \in (0,1/2)$. We suppose that for any $n \geq 1$,
\begin{equation} \label{eligible1}
0 < \eps_n \leq \eps_0.
\end{equation}
Of course, the condition  can be relaxed so that it only holds for    sufficiently large $n$.

Also, fixing $A > 0$, assume that
\begin{equation} \label{eligible2}
\uo \geq -A, \qquad \so^2 \leq A.
\end{equation}
In addition, we assume that for any $n \geq 1$,
\begin{equation} \label{eligible3}
P_{H_n} (u > \uo) = 1, \qquad  P_{H_n} (\sigma^2 \leq A) = 1.
\end{equation}
These conditions are relatively relaxed, except for the second one   in (\ref{eligible3}).   We need this  condition to control the variance of the estimators (whether this
condition can be significantly relaxed is an open question, which we leave  to the future study).    In short,  we focus the study on the class of  marginal densities as follows,
\[
\Lambda_n(\eps_0, A)  = \{f(x) = f(x; \uo, \so, \eps_n, H_n, n) \;  \mbox{has  the form  as in (\ref{f})  that satisfies (\ref{eligible1})-(\ref{eligible3})}     \}.
\]

For any  $t > 0$, it follows  from the triangle inequality that
\[
|\uo(\varphi_n, t)  - \uo| \leq |\uo(\varphi_n, t)  - \uo(\varphi, t)| + |\uo(\varphi, t) - \uo|.
\]
On the right hand side, the first term is the stochastic term, and the second term is  the bias term.  Seemingly,  the performance of the estimator depends on the choice of $t$.  Larger $t$ tends to give a larger stochastic fluctuation but a smaller bias.   It turns out that  the interesting range of $t$ is
$O(\sqrt{\log n})$.  In light of this, we calibrate $t$ through a parameter $\gamma$ by
\[
t = t_n(\gamma)  = \sqrt{    \gamma \log n}, \qquad \gamma  > 0.
\]
We now study the stochastic term and the bias  term separately.
\subsection{The stochastic term}
  We  need the following definition.
\bed
Fixing a constant $r$,    we say that a  sequence $\{b_n\}_{n = 1}^{\infty}$  is $\bar{o}(n^{-r})$ if $n^{r - \delta} |b_n|  \goto 0$ as $n \goto \infty$,  for all $\delta>0$.   Especially, when $r = 0$, we write $\bar{o}(1)$.
\eed

First, we study the stochastic fluctuation of $\varphi_n(t)$ and $\varphi_n'(t)$.
The following lemmas are proved in Section \ref{sec:appen}.
\begin{lemma} \label{lemma:var}
Fix $\eps_0 \in (0,1/2)$,  $A > 0$, and $\gamma \in (0, 1/A)$.  As $n$ tends to $\infty$,
\begin{equation} \label{varbound1}
\sup_{\{f \in \Lambda_n(\eps_0, A)\}}  \{\mathrm{Var}(\varphi_n(t_n(\gamma)))  \}   \leq
 n^{   A \gamma -1},
\end{equation}
and
\begin{equation} \label{varbound2}
 \sup_{\{f \in \Lambda_n(\eps_0, A)\}}  \{\mathrm{Var}(\varphi_n'(t_n(\gamma)))   \}  \lesssim   4 A^2 \gamma  \log(n) \cdot n^{  A \gamma -1}.
\end{equation}
\end{lemma}
The upper bounds in (\ref{varbound1})-(\ref{varbound2})  may be conservative, especially when $\eps_n$ is small.  See the proof for the details (we say two positive sequences $a_n \lesssim b_n$ if $a_n/b_n \leq 1 + o(1)$ for sufficiently large $n$).

We now relate the stochastic fluctuations
of $\uo(\varphi_n; t)$ and $\so^2(\varphi_n; t)$ to  that of  $\varphi_n(t)$ and $\varphi_n'(t)$.  This is achieved by the following lemma.
\begin{lemma}  \label{lemma:fg}
Let $\uo(\cdot; \cdot)$ and $\so^2(\cdot; \cdot)$ be defined as in (\ref{Functional}).  Fix $t > 0$.
For any  differentiable complex-valued  functions  $f$ and $g$ satisfying  $|f(t)| \neq 0$ and $|g(t)| \neq 0$,
\[
|\uo(g, t) - \uo(f, t)|  \leq  \frac{1}{  |f(t)|^2}   \bigl[(    |\uo(g,t)|   \cdot (|f(t) +   |g(t)|)  +\sqrt{2}  |g'(t)|)  |f(t) -g(t)| + \sqrt{2}  |f(t)| \cdot  |(f(t) - g(t))'|,
\]
and
\[
|\so^2(g, t) - \so^2(f, t)|  \leq  \frac{1}{t |f(t)|^2} \cdot \bigl[(|\so^2(g,t)| \cdot t \cdot   (|f(t) + g(t)|) + \sqrt{2} |g'(t)|) \cdot |f(t)  - g(t)|  +  \sqrt{2}  |f(t)| \cdot |(f(t) - g(t))'| \bigr].
\]
\end{lemma}
Apply Lemma \ref{lemma:fg} with $f = \varphi_n$,   $g = \varphi$.       Intuitively,
\[
\varphi_n(t) \approx  \varphi(t),   \qquad  \varphi_n'(t) \approx \varphi'(t).
\]
Also, when $t = t_n(\gamma)$,
\[
\varphi_n(t_n(\gamma)) = \bar{o}(1).
\]
We therefore expect to have
\[
|\uo(\varphi_n, t_n(\gamma)) -  \uo(\varphi; t_n(\gamma))| \leq C \frac{1}{|\varphi(t(\gamma))|} |\varphi_n(t(\gamma)) - \varphi(t(\gamma))| + \frac{|\varphi'(t(\gamma))|}{|\varphi(t(\gamma))|^2} |\varphi_n'(t(\gamma)) - \varphi'(t(\gamma))| \leq \bar{o}(n^{A\gamma -1}).
\]
As a result, we have the following lemma.
\begin{lemma} \label{lemma:stoch}
Fix $\eps_0 \in (0,1/2)$,  $A >  0$,  and   $\gamma   \in (0,1/A)$.  As $n$ tends to $\infty$, except for a probability that tends to $0$,
\[
\sup_{\{f \in \Lambda_n(\eps_0, A)\}}  \{ |\uo(\varphi_n; t_n(\gamma)) - \uo(\varphi; t_n(\gamma))| \}  \leq
\bar{o}(n^{(A\gamma -1)/2}),
\]
and
\[
\sup_{\{f \in \Lambda_n(\eps_0, A)\}}  \{|\so^2(\varphi_n; t_n(\gamma)) -  \so^2(\varphi; t_n(\gamma))|  \}  \leq \bar{o}(n^{(A \gamma -1)/2}).
\]
\end{lemma}
In conclusion, except for a probability that tends to $0$, the stochastic fluctuation of either estimator is of the order  of  $n^{(A \gamma - 1)/2}$. Note that the exponent $(A\gamma -1)/2 < 0$.


\subsection{The bias term}
We now discuss the bias term.   The following lemma is proved in Section \ref{sec:appen}.
\begin{lemma} \label{lemma:r}
Fix $\eps_0 \in (1/2)$ and $A > 0$.    Let $r(t) = r(t; \uo, \so, \eps_n, H_n)$ be as in (\ref{Definer}).  For any $t > 0$ and  $f \in \Lambda_n(\eps_0, A)$, there exists a universal constant  $C > 0$ such that
\begin{align*}
|\uo(\varphi;t)- \uo|  &\leq   C      |r'(t)|,  \\
|\so^2(\varphi; t)  - \so^2| &\leq   C         |r'(t)|/t.
\end{align*}
\end{lemma}

Write for short  $t_n = t_n(\gamma)$.  Under mild conditions,    $r'(t_n) \goto 0$.
We now show some examples where  this is the case.

{\it Example 1. The non-null effects are sparse.}  In this case, we suppose that the parameter $\eps_n$ tends to $0$ as $n$ tends to $\infty$ at a rate faster than that of $1/t_n$.   By the proof of Lemma \ref{lemma:stoch},
\[
|r'(t_n)| \leq \frac{\eps_n}{1 - \eps_n} [\frac{\sqrt{2}}{e t_n} + A t_n].
\]
So as long  as $\eps_n t_n \goto 0$,
$|r'(t_n)| \goto 0$, regardless of the distribution of  $H_n(\cdot, \cdot)$ (of course, the condition of $P_{H_n} (u > \uo) =1$ is still needed).

{\it Example 2.  Elevated means.}   In this case, we suppose that the mean corresponding to the null density is elevated by at least a small amount  $\delta_n > 0$:
\[
P_{H_n}  ( u \geq \uo + \delta_n)  = 1.
\]
Recall that $\uo \geq -A$ and that for any $H_n \in \Lambda_n(\eps_0, A)$,  $P_{H_n} (|\sigma^2 - \so^2| \leq A)  = 1$.
Similar to the proof  of Lemma \ref{lemma:stoch},
\[
|r'(t_n)| \leq    \frac{\eps_n}{1 - \eps_n}     (\delta_n + \frac{\sqrt{2}}{e t_n}  + A t_n) e^{-\delta_n  t_n /\sqrt{2}}.
\]
As a result, we have the following lemma, whose proof is elementary so we omit it.
\begin{lemma} \label{lemma:exp1}
If there is some constant $c_0 > 0$ such that
\begin{equation} \label{Definedelta}
\liminf_{n \goto \infty}   \bigl\{  \frac{ \delta_n t_n}{ \sqrt{2}  \log(t_n) }  \bigr\}   \geq  (c_0 + 1),
\end{equation}
then $|r'(t_n)|  \lesssim  A \eps_n t_n^{-c_0}$.
\end{lemma}
As a result,   as $n \goto \infty$,   the bias $\goto  0$ if (\ref{Definedelta}) holds for some constant  $c_0 > 0$,  whether
$\eps_n$   tends to $0$ or not.

{\it Example 3.  When the bivariate random variables $(u, \sigma^2)$ have a smooth joint density.}    We re-center $u$ and $\sigma^2$ by letting
$\delta = u - \uo$ and $\kappa = (\sigma^2 - \so^2)/2$.  Denote the joint density of $(\delta, \kappa)$ by $h_n(\cdot, \cdot)$.  We show that the  $r'(t_n) = o(1)$  under mild   smoothness conditions on  $h_n(\cdot, \cdot)$.
In detail,   for each fixed $\delta > 0$,  let $h_n^{FT}(\cdot | \delta)$ be the Fourier transform of  the conditional density  $h_n(\kappa | \delta)$.  Fix $\alpha > 0$.  Suppose that  there is a generic  constant $C >  0$ such that
for all $\delta$ in the range,
\begin{equation} \label{FT}
|h_n^{FT}(t_n|\delta)| \leq C (1 + |t_n|)^{-\alpha}, \qquad |\frac{d}{dt} h_n^{FT}(t_n|\delta)| \leq C(1 +  |t_n|)^{-(\alpha  + 1)}.
\end{equation}
We have the following lemma, whose proof is elementary  so we omit it.
\begin{lemma} \label{lemma:exp2}
Suppose (\ref{FT})  holds for some constant $C > 0$ and $\alpha > 0$. Then there is a generic constant $C > 0$ such that
\[
|r'(t_n)|  \leq   C  \eps_n  |t_n|^{-2(\alpha + 1)},
\]
\end{lemma}
Note that $r'(t_n) \goto 0$ in a much  broader setting than that in this example.

Combining Lemma \ref{lemma:stoch}-\ref{lemma:r} and the above examples, we have the following theorem.
\begin{thm} \label{thm:Main}
Fix $\eps_0 \in (0,1/2)$,  $A >  0$,  and   $\gamma   \in (0,1/A)$.
Suppose that when $n$ tends to $\infty$,  at least one of the three conditions below holds:
\begin{description}
\item[a.]  $\lim_{n \goto \infty} (\eps_n \cdot t_n(\gamma)) = 0$,
\item[b.]   $P_{H_n} (u > \uo + \delta_n) =1$, where $\delta_n$ satisfies (\ref{Definedelta}) for some constant $c_0 > 0$.
\item[c.]   (\ref{FT}) holds for some parameter $\alpha > 0$.
\end{description}
Then the estimators $\uo(\varphi_n; t_n(\gamma))$ and $\so^2(\varphi_n; t_n(\gamma))$ are
consistent with respect to the null parameters $\uo$ and $\so^2$, respectively, uniformly across all densities in $\Lambda_n(\eps_0, A)$ that  satisfy one or more of the conditions (a), (b), and (c).
\end{thm}
We remark that while choosing $\gamma \in (0,1/A)$ ensures consistency,
 different choices of $\gamma$ affect the convergence rate of the estimators.  The optimal choice of $\gamma$ depends on unknown parameters and is hard to set.   In Section \ref{sec:simul}, we investigate how to choose $\gamma$ with simulated data. In our experience,  when $A$ is not very large, it is usually appropriate to choose $\gamma \approx 0.2$.

We  also remark that in Theorem \ref{thm:Main} (as well as Theorems \ref{thm:eps1}, \ref{thm:eps2} below), we have assumed independence of the test statistics $X_j$.  When the test statistics $X_j$ are correlated, the bias  of the estimators remain the same,  but the variance of the estimators may inflate by a factor. On the other hand, if the correlation is relatively weak, the estimators continue to perform well. In Section \ref{sec:simul}, we investigate an simulation example with block-wise
dependence among $X_j$. The simulation results suggest that the estimators continue to  perform  well
when the block size is small (e.g. $\leq 100$). See Section \ref{sec:simul} and Figure \ref{Fig:dependent} for the details.

\section{Estimating the proportion of non-null effects}   \label{sec:prop}

The proportion has an identifiability issue that is very similar to
that of the null parameters.  The issue can also be resolved
similarly in GEV and GEM.  In \cite{JC, Jin, CJ}, we have carefully
investigated the problem of estimating the proportion in GEV.
Similarly to estimating the null parameters,  a Fourier approach was
introduced (e.g. \cite{Jin, JC}). Compared to existing approaches in
the literature  (e.g. \cite{Efronetal, Langaas, Rice, Storey1,
Swanepoel, Wasserman}),  the Fourier approach was proven to be
successful in a much broader setting. Especially,  it  was shown to
be successful  without the so-called {\it purity} condition, a
notion introduced in \cite{Wasserman}. Later in \cite{CJ}, the
approach was shown to also attain the optimal rate of convergence
over a wide class of situations.

We now shift our attention to GEM.  Despite the encouraging development,  the Fourier approach in \cite{Jin, JC}  ceases to perform well in this case.     In fact,    in this case,    it can be shown that none of the aforementioned approaches is uniformly consistent with the proportion.
Therefore, it is necessary to develop a new approach.

In this section, we propose a new approach to estimating the proportion by using the generalized Fourier transformation, as a natural extension of  the ideas in preceding sections.  We discuss two cases separately:  the case where  the null parameters are known, and the case where  the null parameters are unknown. In both cases,
we show that under mild conditions, the proposed  approach  is   uniformly  consistent with the proportion.

\subsection{Known null parameters}
Recall that
\[
\varphi_0(t) = (1 - \eps_n) e^{ \omega \uo t  + i \so^2 t^2/2}.
\]
The key observation is that, when the null parameters $(\uo, \so^2)$ are known,  $\eps_n$ can be easily solved from $\varphi_0(t)$ by
\[
\eps_n  \equiv    1 -  e^{- \omega \uo t  - i \so^2 t^2/2}  \varphi_0(t).
\]

Inspired by this, we introduce the functional
\begin{equation}  \label{Defineeps}
\eps_n(g; t,  \uo, \so^2)   =   1 -   e^{- \omega \uo t  - i \so^2 t^2/2}  g(t).
\end{equation}
where $g(t)$ is any complex-valued function.    Recall that
\[
\varphi_n(t) \approx  \varphi(t) \approx \varphi_0(t).
\]
By the continuity of the functional, we hope that for an  appropriately chosen $t$,
\[
\eps_n(\varphi_n; t, \uo, \so^2) \approx  \eps_n(\varphi_0; t, \uo, \so^2)  \equiv \eps_n.
\]

We now analyze the variance and the bias of this estimator.   As before, let $t_n(\gamma) = \sqrt{  \gamma \log n}$.   For any $H_n(\cdot, \cdot)$ satisfying $P_{H_n}(\sigma^2 \leq A) = 1$,  direct calculations show that
\[
\mathrm{Var}( \eps_n(\varphi_n; t_n(\gamma), \uo, \so^2) \leq  \frac{1}{n} E[e^{ -\sqrt{2} t_n(\gamma) (X_1 - \uo)}] \leq   \frac{1}{n}  [(1 - \eps_n) n^{\so^2 \gamma}  + \eps_n n^{ A  \gamma}],
\]
so the standard deviation of the estimator is of the order  of $o(\eps_n)$ when
\[
n^{A\gamma -1}  = o(\eps_n^2).
\]
At the same time, by elementary calculus,   the bias of the estimator equals to
\[
\bigl|E[ \eps_n(\varphi_n; t_n, \uo, \so^2)]  - \eps_n \bigr|  =  \eps_n \cdot  \bigl|  \int e^{\omega (u - \uo) t_n + i (\sigma^2 - \so^2) t_n^2/2} d H_n(u, \sigma) \bigr|,
\]
which is of the order of $o(\eps_n)$ if either of the aforementioned  conditions  (b)  or (c) holds.  Combining these gives the following theorem.
\begin{thm}  \label{thm:eps1}
Fix  $\uo$,  $A > 0$,  $\so^2  \in (0, A)$,  $\gamma  \in (0, 1/A)$, and a sequence of positive numbers $b_n$ satisfying $\lim_{n \goto \infty} b_n  = 0$.
Consider a sequence of parameters  $\eps_n  \in (0,1)$ and a sequence of  bivariate distribution $H_n(u, \sigma)$ such that for sufficiently large $n$,   $P_{H_n}( u > \uo, \sigma^2 \leq A) = 1$ and  $\eps_n^{-2} n^{A \gamma - 1}   \leq b_n$.
Also, suppose that when $n$ tends to $\infty$,  at least one of the two conditions below holds:
\begin{description}
\item[b.]   $P_{H_n} (u > \uo + \delta_n) =1$, where $\delta_n$ satisfies (\ref{Definedelta}) for some constant $c_0 > 1$.
\item[c.]   (\ref{FT}) holds for some parameter $\alpha  > 0$.
\end{description}
Then as $n \goto \infty$,   except for a probability that  tends to zero,
\[
\bigl|\frac{\eps_n(\varphi_n; t_n(\gamma), \uo, \so^2)}{\eps_n} - 1 \bigr| \goto 0,
\]
uniformly for all $\eps_n$ and $H_n(\cdot, \cdot)$ satisfying the conditions  above.
\end{thm}
In other words,  $\eps_n(\varphi_n; t_n(\gamma))$  is uniformly consistent with $\eps_n$ provided that either (b) or (c) holds, and that the variance of the estimator is of a  smaller order than that of  $\eps_n^2$.
The latter is satisfied when $\eps_n$ tends to $0$ slowly enough.

\subsection{Unknown null parameters}
When the null parameters are unknown,  a natural approach is to  estimate the null parameters using the approach in Section \ref{sec:main} first, then plug in the estimated values to estimate the
proportion. In other words,  we first estimate the null parameters by
\begin{equation} \label{ushat}
  \hat{u}_0(\gamma) =  \uo(\varphi_n; t_n(\gamma)), \qquad \hat{\sigma}_0^2(\gamma)  = \so^2(\varphi_n; t_n(\gamma)).
\end{equation}
We then estimate the proportion by the plugging estimator,
\[
\eps_n(\varphi_n; t_n(\gamma), \hat{u}_0(\gamma), \hat{\sigma}_0^2(\gamma)) =   1  -  e^{ - \omega \hat{u}_0(\gamma)
t_n(\gamma) - i \hat{\sigma}_0^2(\gamma)  t_n^2(\gamma)/2}
\varphi_n(t_n(\gamma)).
\]
Note that the bias  of  both $\hat{u}_0(\gamma)$ and $\hat{\sigma}_0^2(\gamma)$ are typically of the order of  $o(\eps_n)$, and their  variance   are of the same order as that of  $\eps_n(\varphi_n; \gamma)$. Therefore, replacing $(\uo, \so^2)$  by $(\hat{u}_0(\gamma), \hat{\sigma}_0^2(\gamma))$ does not increase either the
bias or variability of the estimator.    The following
theorem is proved in Section \ref{sec:appen}.
\begin{thm}  \label{thm:eps2}
Fix  $\uo$,  $A > 0$,  $\so^2  \in (0, A)$,  $\gamma  \in (0, 1/A)$, and a sequence of positive numbers $b_n$ satisfying $\lim_{n \goto \infty} b_n  = 0$.
Consider a sequence of parameters  $\eps_n  \in (0,1)$ and a sequence of  bivariate distribution $H_n(u, \sigma)$ such that for sufficiently large $n$,   $P_{H_n}( u > \uo, \sigma^2 \leq A) = 1$ and  $\eps_n^{-2} n^{A \gamma - 1}   \leq b_n$.
Also, suppose that when $n$ tends to $\infty$,  at least one of the two conditions below holds:
\begin{description}
\item[b.]   $P_{H_n} (u > \uo + \delta_n) =1$, where $\delta_n$ satisfies (\ref{Definedelta}) for some constant $c_0 > 1$.
\item[c.]   (\ref{FT}) holds for some parameter $\alpha  > 0$.
\end{description}
Then as $n \goto \infty$,   except for a probability that  tends to zero,
\[
\bigl|\frac{ \eps_n(\varphi_n; t_n(\gamma), \hat{u}_0(\gamma),  \hat{\sigma}_0^2(\gamma))}{\eps_n} - 1 \bigr| \goto 0,
\]
uniformly for all $\eps_n$ and $H_n(\cdot, \cdot)$ satisfying the conditions  above.
\end{thm}

\section{Simulations} \label{sec:simul}
 In this section, we conduct simulation studies for investigating the performance of the proposed estimators of $(\uo, \so^2, \eps)$ with a finite $n$. We write for short
\[
\hat{u}_0(\gamma) = \uo(\varphi_n; t_n(\gamma)),  \qquad \hat{\sigma}_0^2(\gamma)   = \so^2(\varphi_n; t_n(\gamma)), \qquad  \hat{\eps}_n(\gamma) =  \eps_n(\varphi_n; t_n(\gamma), \hat{u}_0(\gamma), \hat{\sigma}_0^2(\gamma)).
\]
  Specifically, we are interested
in  four aspects: (1) how different choices of $\gamma$ affect the
estimation errors  of $\hat{u}_0(\gamma)$,
$\hat{\sigma}_0^2(\gamma)$ and $\hat{\eps}_n(\gamma)$; and what
$\gamma$ values we should recommend in practice; (2) the effect of
different choices of the proportion $\eps$ and the mixing
distribution $H_n(\cdot, \cdot)$; (3)  the effect of   larger $n$; and
(4) the effect of dependent structures.\vspace{5pt}

\begin{figure}[htb]
\begin{center}
\includegraphics[height=5.8in, angle=270]{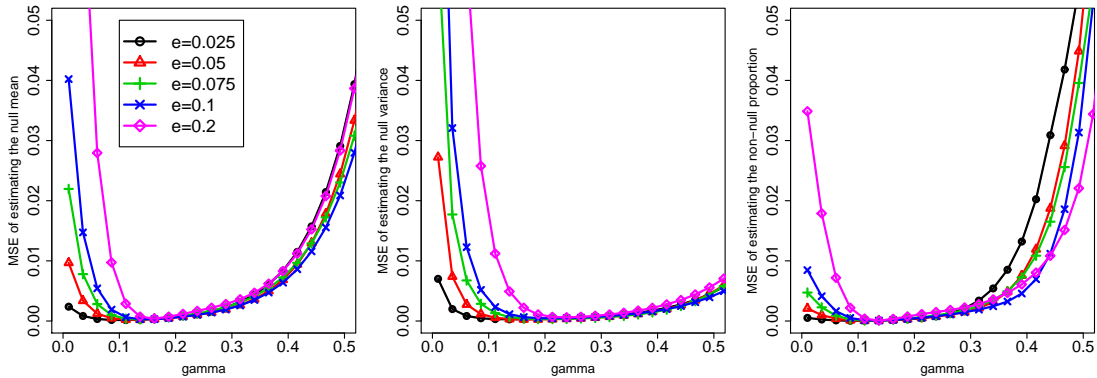}
\end{center}
\caption{MSE for $\hat{u}_0(\gamma)$ (\textbf{left}),
$\hat{\sigma}_0^2(\gamma)$ (\textbf{middle}) and
$\hat{\eps}_n(\gamma)$ (\textbf{right})
for $100$ repetitions. The $x$-axis displays $\gamma$,  and the
$y$-axis displays the MSE. The null parameters $(\uo, \so^2) = (-1,
1)$.
Different colors of the curves represent different values of
$\eps$.}
\label{Fig:varyingEps}
\end{figure}

{\it \underbar{Example 1.}  Different choices of the tuning
parameter $\gamma$}.   In this example, we let $n=50,000$,  $(\uo,
\so^2) = (-1, 1)$,   and $\eps = 0.025 \times (1, 2, 3, 4, 8)$. We
choose $20$ different $\gamma$ ranging  from $0.01$ to $0.5$ with
equal inter-distances. For each combination of $(\eps, \gamma)$, we
conduct an experiment with the following four steps.
\begin{itemize}
\item {\it Step 1.}    For each $1 \leq j \leq n(1 - \eps)$, draw $X_j  \sim N(u_0,
\sigma_0^2)$ to represent a null effect.
\item {\it Step 2.}   For each $n(1 - \eps) + 1 \leq j \leq n$,   draw independently a sample $u \sim \mathrm{Uniform}(1, 2)$ and a
sample $\sigma \sim \mathrm{Uniform}(0.5, 1.5)$.  Then, draw $X_j
\sim N(u, \sigma^2)$ to represent a non-null effect.
\item {\it Step 3.}  Calculate  $\hat{u}_0(\gamma), \hat{\sigma}_0^2(\gamma)$, and $\hat{\eps}_n(\gamma)$.   \item {\it Step 4.}  Repeat  Steps 1--3  for 100 times.
\end{itemize}
The results are reported in Figure \ref{Fig:varyingEps}, from which
we can see that the MSE are the smallest when  $\gamma\in (0.15,
0.25)$. Also, the MSE are not  sensitive to different choices of
$\gamma$: they remain about the same for different $\gamma \in
(0.15, 0.25)$. All of the three estimators     $\hat{u}_0(\gamma),
\hat{\sigma}_0^2(\gamma)$, and $\hat{\eps}_n(\gamma)$ have
satisfactory performances: when $\gamma=0.2$, the  MSE are as small
as the order of $10^{-4}$. Somewhat surprisingly, in this example,
different $\eps$ do not have a prominent effect on the MSE.
\vspace{5pt}

{\it \underbar{Example 2.}  The effect of different mixing
distribution $H_n(\cdot, \cdot)$.}   In this example, we set $n =
50,000$, $(\uo, \so^2, \eps) = (-1, 1, 0.05)$, and choose $20$
different $\gamma$ ranging from $0.01$ to $0.5$ with equal
inter-distances. Compared to Example 1, we conduct experiments with
different choices of the mixing distribution $H_n(\cdot, \cdot)$.   We
consider two scenarios. In the first scenario,  independently,   $(u
- \uo)  \sim \mathrm{Gamma}(10, 0.25)$  ($\mathrm{Gamma}(k, \theta)$
is the Gamma distribution with shape parameter $k$ and scale
parameter $\theta$), and  $\sigma  \sim \mathrm{Uniform}(0.5, 1.5)$.
The parameters $(10, 0.25)$ are chosen such that  the mean value of
the random variable $u$ is $1.5$,   the same as that  in the
preceding example. In the second scenario,   independently,     $u
\sim  \mathrm{Uniform}(1, 2)$ and $\sigma  \sim \mathrm{Gamma}(10,
0.1)$.

\begin{figure}[htb]
\begin{center}
\includegraphics[height=2.5in, angle=270]{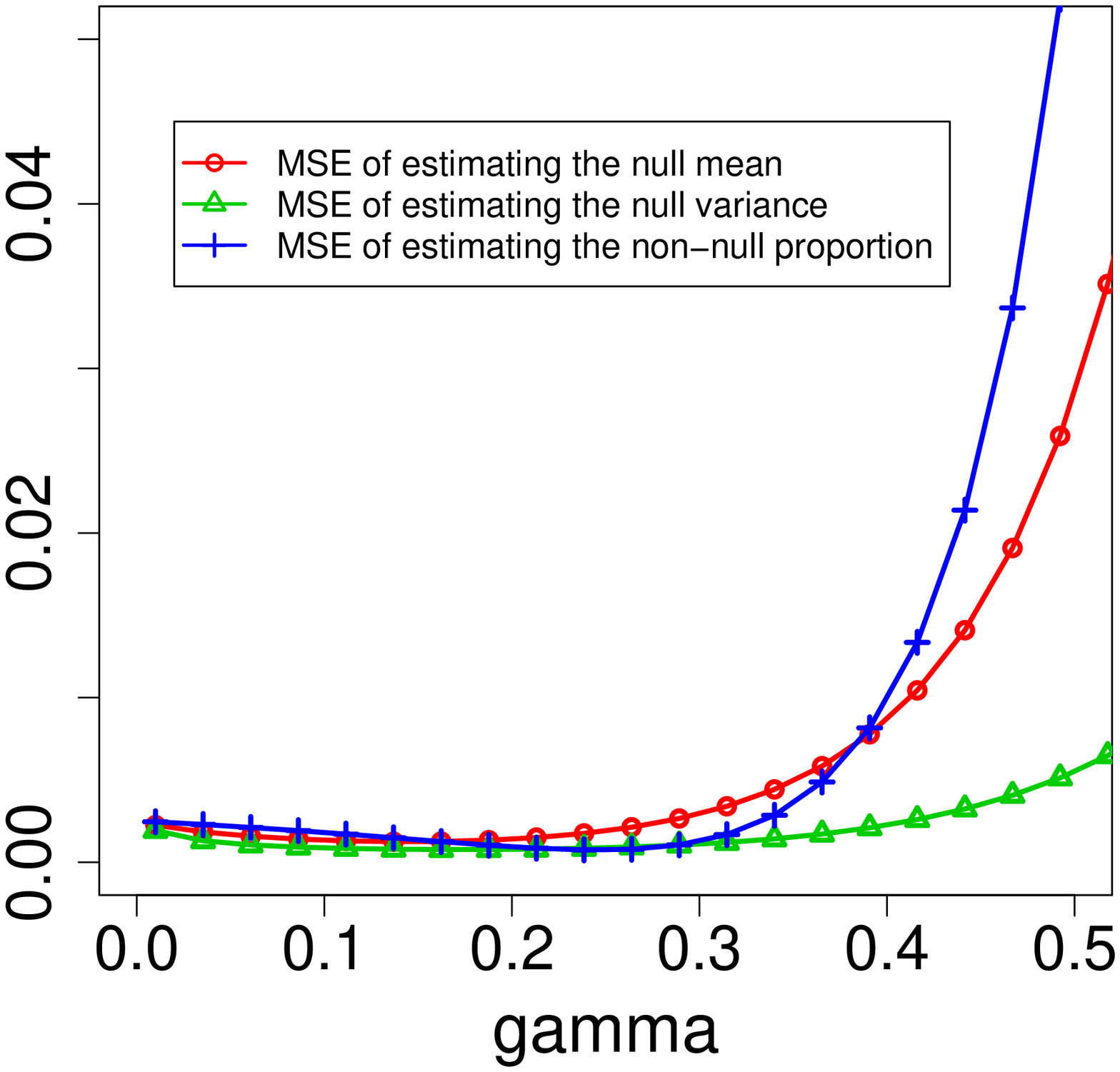}
\includegraphics[height=2.5in, angle=270]{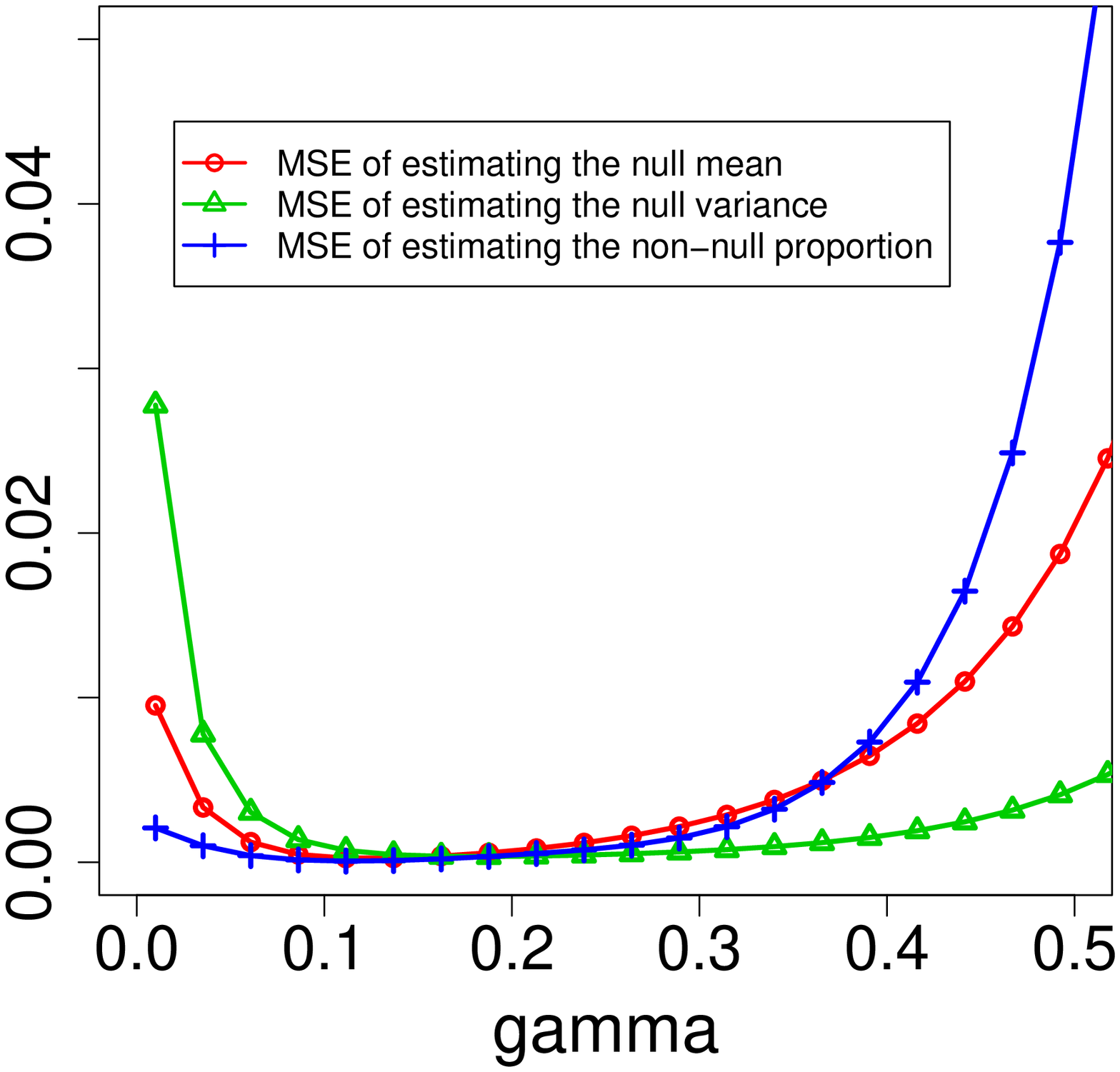}
\caption{MSE for $\hat{u}_0(\gamma)$ (\textit{red}),
$\hat{\sigma}_0^2(\gamma)$ (\textit{green}) and
$\hat{\eps}_n(\gamma)$ (\textit{blue})
for Scenario 1 (\textbf{left}) and Scenario 2 (\textbf{right})
considered in Example 2. The $x$-axis displays $\gamma$ and the
$y$-axis displays the MSE. $n = 50,000$ and $(\uo, \so^2, \eps) =
(-1, 1, 0.05)$.}\label{Fig:alter}
\end{center}

\end{figure}

For each scenario and each $\gamma$,   we  run experiments following
Steps 1--4 as in Example 1,  but  with the current choice of
$H_n(\cdot, \cdot)$. The MSE for  $\hat{u}_0(\gamma),
\hat{\sigma}_0^2(\gamma)$ and $\hat{\eps}_n(\gamma)$ are reported in
Figure \ref{Fig:alter}.   From this figure,  a similar conclusion
can be drawn: the estimators perform well in both scenarios, with
the MSE  as small as   $10^{-4}$--$10^{-3}$. The best range of
$\gamma$ is $(0.15, 0.2)$. In this range, the MSE is relatively
insensitive to different choices of $\gamma$.

 It is noteworthy that in the first scenario,    the support of
random variable  $u$ is not bounded away from the null parameter
$\uo$.  It is also noteworthy that in the second scenario, $\sigma$
is unbounded such that the assumption (\ref{eligible3}) is violated.
Despite the seeming challenges in these two scenarios,  the proposed
approach continues to perform well. This suggests that the proposed
approaches are successful in a broader situations than that
considered in Sections \ref{sec:main} and \ref{sec:prop}.\vspace{5pt}

{\it \underbar{Example 3}.  The effect of larger $n$.}   In this
example, we fix $(\uo, \so^2, \eps) = (-1, 1, 0.05)$. Since the MSE
is relatively insensitive to different choices of $\gamma$, we fix
$\gamma = 0.2$.  For the mixing distribution $H_n(\cdot, \cdot)$,  we
let $u \sim \mathrm{Uniform}(1, 2)$ and $\sigma \sim
\mathrm{Uniform}(0.5, 1.5)$, independently of each other. According
to the asymptotic analysis in preceding sections,  we understand
that the performance of proposed estimators improves when $n$
increases.  In this example, we validate this point by choosing
$n=10^4 \times (1, 3, 5, 8, 10)$.   For each $n$,  we run
experiments following Steps 1--4 as in Example 1.    The results are
summarized in Table~\ref{Table:varyingN}.   The MSE of all
$\hat{u}_0(\gamma)$, $\hat{\sigma}_0^2(\gamma)$ and
$\hat{\eps}_n(\gamma)$ decreases as $n$ increases.  This fits well
with the asymptotic analysis in Sections
\ref{sec:main} and \ref{sec:prop}.\vspace{5pt}

\begin{table}
\begin{center}
\caption{MSE for $\hat{u}_0(\gamma)$, $\hat{\sigma}_0^2(\gamma)$, and $\hat{\eps}_n(\gamma)$
for different $n$, where we take $\gamma = 0.2$. The parameters   $(\uo,
\so^2, \eps) = (-1, 1, 0.05)$.  For the mixing distribution $H_n(u,
\sigma)$,   $u \sim \mathrm{Uniform}(1,2)$ and $\sigma \sim
\mathrm{Uniform}(0.5, 1.5)$ independently. In each cell, the MSE
equals  the cell value  times $10^{-4}$. }
\begin{tabular}{c|ccccc}\hline
n & $10^4$ & $3\times 10^4$ & $5\times 10^4$ & $8\times 10^4$
&$10^5$ \\\hline MSE for $\hat{u}_0(\gamma)$  & 41.28 & 10.46& 5.66 &
4.01 & 2.73 \\
MSE for $\hat{\sigma}_0^2(\gamma)$  & 16.47 & 6.93 & 2.36 & 1.81 & 1.48\\
MSE for $\hat{\eps}_n(\gamma)$ & 20.13 & 5.28 & 4.17 & 2.87 & 2.01
\\\hline
\end{tabular}
\label{Table:varyingN}
\end{center}
\end{table}

\begin{figure}[htb]
\begin{center}
\includegraphics[height=2.5in, angle=270]{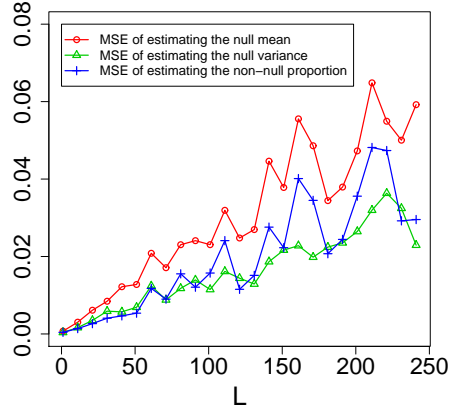}
\caption{MSE for $\hat{u}_0(\gamma)$ (\textit{red}),
$\hat{\sigma}_0^2(\gamma)$ (\textit{green}) and
$\hat{\eps}_n(\gamma)$ (\textit{blue})
when the data $X_j$ are block-wise dependent with a block size $L$
(displayed in the $x$-axis; see Example 4 for the details).  The
parameters $(\uo, \so^2, \eps) = (-1, 1, 0.05)$.  For the mixing
distribution $H_n(u, \sigma)$,   $u \sim \mathrm{Uniform}(1,2)$ and
$\sigma \sim \mathrm{Uniform}(0.5, 1.5)$ independently. }
\label{Fig:dependent}
\end{center}
\end{figure}

{\it \underbar{Example 4}.  The effect of dependence.}   In this
example, we fix  $n = 50,000$, $(\uo, \so^2, \eps) = (-1, 1, 0.05)$.
We investigate how the dependent structures may affect the
performance of the proposed procedures.  For each $L$ ranging  from
$1$ to $250$ with an increments of $10$,  we generate samples as
follows.
\begin{enumerate}
\item For each $1 \leq j \leq n(1 - \eps)$, set $(\mu_j, \sigma_j) = (\uo, \so)$.
\item For each $n(1 - \eps) + 1 \leq j \leq n$,  draw $\mu_j \sim \mathrm{Uniform}(1,2)$ and $\sigma_j \sim \mathrm{Uniform}(0.5, 1.5)$.
\item Draw $w_1,  \cdots, w_{n+L}$ independently from $N(0,1)$. For $1 \leq j \leq n$,   let $z_j= \sum_{k=j}^{k=j+L} \frac{w_k}{\sqrt{L+1}}$.
Note that marginally $z_j\sim N(0,1)$.
\item For $1 \leq j \leq n$, let $X_j=\mu_j+\sigma_j\cdot z_j$.
\end{enumerate}
The data generated in this way is block-wise dependent, with the
block size being controlled by $L$.   Fix $\gamma = 0.2$. We
calculate $\hat{u}_0(\gamma)$, $\hat{\sigma}_0^2(\gamma)$, and
$\hat{\eps}_n(\gamma)$,   and repeat the experiment for $100$ times.
We then calculate the MSE. The results are summarize in Figure
\ref{Fig:dependent}. While the MSE increase with the block size $L$,
we also note that the MSE remain small when, say,
 $L\leq 50$ (all three curves fall below $0.02$).  This suggests that  the proposed
methods are relatively robust for short-range dependence.

\section{Proofs} \label{sec:appen}
\subsection{Proof of Lemma \ref{lemma:GEV} and \ref{lemma:GEM} }
We prove Lemma \ref{lemma:GEV} first. Consider two density functions $f_k(x) =  f_k(x; \uo^{(k)}, (\so^2)^{(k)}, \eps^{(k)}, H^{(k)})$ that satisfy (\ref{Definef})-(\ref{GEV}), $k = 1, 2$.   For short, denote $(u_k, \sigma_k^2, \eps_k, H_k)  =  (\uo^{(k)}, (\so^2)^{(k)}, \eps^{(k)}, H^{(k)})$.    Suppose $f_1 = f_2$.       We want to show that
$(u_1,\sigma_1,\epsilon_1)=(u_2,\sigma_2,\epsilon_2)$.   Note that the Fourier transformation of $f_1$ and $f_2$ must be identical. By direct calculations,  with $s_i(t)$ as defined in (\ref{Definess}),
\begin{equation} \label{GEVpf0}
(1-\epsilon_1)e^{itu_1-\frac{\sigma_1^2t^2}{2}}(1+s_1(t))=(1-\epsilon_2)e^{itu_2-\frac{\sigma_2^2t^2}{2}}(1+s_2(t)).
\end{equation}
We first show $\sigma_1 = \sigma_2$.     By (\ref{GEVpf0}),
\begin{equation} \label{GEVpf1}
e^{- \frac{(\sigma_2^2 - \sigma_1^2) t^2}{2}}   = \biggl|  e^{it (u_1 - u_2)}  \frac{(1 - \eps_1) (1 + s_1(t))}{(1 - \eps_2) (1 + s_2(t))} \biggr|.
\end{equation}
Note that $|s_k(t)| \leq \epsilon_k/(1-\epsilon_k) < 1$, where  $|\eps_k| \leq \eps_0 < 1/2$.  Therefore,   the right hand side of (\ref{GEVpf1}) is bounded away from both $0$ and $\infty$ by a constant. Letting $t$ tend to $\infty$ implies that $\sigma_1 = \sigma_2$.

Next, we show $(u_1, \eps_1) = (u_2, \eps_2)$.  By (\ref{GEVpf0}) and $\sigma_1 = \sigma_2$,
\begin{equation} \label{GEVpf2}
(1-\epsilon_1) (1+s_1(t))=(1-\epsilon_2)e^{it (u_2 - u_1)}(1+s_2(t)).
\end{equation}
Fix a small positive number $a > 0$, let $\phi_a(t)$ be the density of $N(0, 1/a)$.  Times $\phi_a(t)$ to both sides of (\ref{GEVpf2}) and integrate in terms of $t$.  By direct calculations and Fubini's theorem,    the left hand side of (\ref{GEVpf2}) is
\begin{equation} \label{GEVpf3}
(1-\epsilon_1)+\epsilon_1\int \frac{a}{\sqrt{\sigma^2-\sigma_1^2+a^2}}  \mathrm{exp} \bigl(-\frac{(u-u_1)^2}{2(\sigma^2-\sigma_1^2+a^2)} \bigr)d H_1(u,\sigma),
\end{equation}
and the right hand side of (\ref{GEVpf2}) is
\begin{equation} \label{GEVpf4}
(1-\epsilon_2) e^{-\frac{(u_2-u_1)^2}{2a^2}}+\epsilon_2\int \frac{a}{\sqrt{\sigma^2-\sigma_2^2+a^2}} \mathrm{exp} \bigl(- \frac{(u-u_1)^2}{2(\sigma^2-\sigma_2^2+a^2)} \bigr)d H_2(u,\sigma).
\end{equation}
Note that by  Dominant Convergence  Theorem  (DCT),  any fixed $H(\cdot, \cdot)$ satisfying $P_H(\sigma \geq \so) = 1$ and $P_{H} \bigl( (u, \sigma)  =  (\uo, \so) \bigr)   = 1$,
\begin{equation} \label{GEVpf5}
\lim_{a \goto 0}   \int \frac{a}{\sqrt{\sigma^2-\sigma_0^2+a^2}}  \mathrm{exp}( -\frac{(u-u_0)^2}{2(\sigma^2-\sigma_0^2+a^2)}) d H(u,\sigma)  = 0.
\end{equation}
Combining (\ref{GEVpf3})--(\ref{GEVpf5})  gives  $(1 - \eps_1) = \lim_{a \goto \infty} [ (1 - \eps_2) \mathrm{exp}(- \frac{(u_2 - u_1)^2}{2a^2})]$,  which immediately implies $(u_1, \eps_1) = (u_2, \eps_2)$.  This proves Lemma \ref{lemma:GEV}.

Consider Lemma \ref{lemma:GEM}.  The difference is now that both $f_i$ satisfy (\ref{Definef})-(\ref{GEM}).
Similarly, suppose  $f_1 \equiv f_2$.   We want to show that
$(u_1,\sigma_1,\epsilon_1)=(u_2,\sigma_2,\epsilon_2)$.    By direct calculations,
the generalized Fourier transform of $f_k$ are $(1 - \eps_k) e^{\omega u_k  t  + i \sigma_k^2 t^2}[ 1  + r_k(t)]$, $k = 1, 2$. It follows that
\begin{equation} \label{GEVpf6}
 e^{\omega  t (u_1 - u_2)}  \cdot  e^{ i (\sigma_1^2  - \sigma_2^2) t^2}  =   \biggl( \frac{1 - \eps_2}{1 - \eps_1}  \biggr)\biggl(  \frac{1 + r_2(t)}{1 + r_1(t)} \biggr).
\end{equation}
Let $t \goto \infty$ on both sides.
By the condition of $P_{H_k} (u > u_k) = 1$,    $r_k(t) \goto 0$. Comparing the modules of both sides gives      $u_1 = u_2$ and $\eps_1 = \eps_2$.  Combining this with (\ref{GEVpf6}),    $e^{ i (\sigma_1^2  - \sigma_2^2) t^2}  =     \frac{1 + r_2(t)}{1 + r_1(t)}$.  Letting $t \goto \infty$, the right hand side tends to $1$.   Therefore,  $\sigma_1 = \sigma_2$.  \qed

\subsection{Proof of Lemma \ref{lemma:var}}
It is sufficient to show  that for any $t > 0$ and $f \in \Lambda_n(\eps_0, A)$,
\begin{equation} \label{varpf1}
 \mathrm{Var}(\varphi_n(t))   \leq
\frac{1}{n} e^{- \sqrt{2} \uo t + \so^2 t^2}  [(1 - \eps_n) + \eps_n e^{(A - \so^2) t^2}],
\end{equation}
and
\begin{equation} \label{varpf2}
\mathrm{Var}(\varphi_n'(t)) \lesssim   \frac{1}{n} e^{ - \sqrt{2} \uo t   +
\so^2 t^2} [(1 - \eps_n)   (\so^2 +2 \uo^2  + 4\so^4  t^2)  + \eps_n (A + 2 \uo^2 + 4 A^2 t^2)
e^{ (A - \so^2) t^2}].
\end{equation}
In fact, once these are proved, the claim follows from  (2.12)-(2.14) by taking $t = t_n(\gamma)$.

Consider (\ref{varpf1}).
Direct calculations show that
\[
\mathrm{Var}(\varphi_n(t))  \leq  \frac{1}{n} E[ |e^{\omega t X_j}|^2]   \leq
\frac{1}{n}  E[e^{- \sqrt{2} t X_j}].
\]
Direct calculations show that
\[
E[e^{- \sqrt{2} t X_j}]   = e^{- \sqrt{2} t \uo + \so^2 t^2} [(1 - \eps_n) + \eps_n \int e^{-\sqrt{2} (u - \uo) t + (\sigma^2 - \so^2) t^2} d H_n(u, \sigma)].
\]
The claim follows by $u_0 >  -A$,  $\so^2 \leq A$,  and $P_{H_n} (u > \uo,  \so^2 \leq A) = 1$.

Consider (\ref{varpf2}). Similarly,
\[
\mathrm{Var}(\varphi_n'(t)) \leq \frac{1}{n} E[  | \omega X_j e^{ \omega t X_j}|^2] \leq \frac{1}{n} E [X_j^2 e^{- \sqrt{2} t X_j}].
\]
By direct calculations,
\begin{equation} \label{var01}
E[X_j^2 e^{ - \sqrt{2} t X_j}] =  I + II,
\end{equation}
where
\[
I =  (1 - \eps_n) [\so^2 + (- \uo + \sqrt{2} \so^2 t)^2]  e^{-\sqrt{2} \uo t + \so^2 t^2},
\]
and
\[
II = \eps_n \int [\sigma^2 + (- u + \sqrt{2} \sigma^2 t)^2] e^{-\sqrt{2} u t + \sigma^2 t^2} d H_n(u, \sigma).
\]
By Schwartz inequality,
\[
(- \uo + \sqrt{2} \so^2 t)^2 \leq  2(\uo^2 + 2 \so^4 t^2),
\]
\[
(-u  +  \sqrt{2} \sigma^2 t^2)^2   = (- \uo - (u - \uo)  + \sqrt{2} \sigma^2 t^2)^2 \leq  2(\uo^2  +  2 \sigma^2 t^2 + (u - \uo)^2).
\]
So
\begin{equation} \label{var02}
I \leq  (1 - \eps_n) [\so^2 + 2 \uo^2 +  4 \sigma_0^4 t^2] e^{-\sqrt{2} \uo t + \so^2 t^2},
\end{equation}
and
\begin{equation} \label{var1}
II \leq   \eps_n [ \int (\sigma^2 + 2 \uo^2 + 4 \sigma^4 t^2) e^{ -\sqrt{2}  u  t + \sigma^2 t^2}  d H_n(u, \sigma) +  2   e^{-\sqrt{2} \uo t}  \int (u - u_0)^2 e^{-\sqrt{2} (u - \uo) t + \sigma^2 t^2} d H_n(u, \sigma)].
\end{equation}
Note that    $\sup_{x > 0} 2x^2 e^{-\sqrt{2} t x}  = 2/(et^2)$.  It follows
\begin{equation} \label{var2}
\int (u - u_0)^2 e^{-\sqrt{2} (u - \uo) t + \sigma^2 t^2} d H_n(u, \sigma)]  \leq [2/(e t^2)] \int  e^{\sigma^2 t^2} d H_n(u, \sigma).
\end{equation}
Inserting   (\ref{var2})  into (\ref{var1})  and recalling that $P_{H_n} (u > \uo,  \sigma^2 \leq A) = 1$,
\begin{equation} \label{var3}
II \leq \eps_n  [A + 2 \uo^2 + 4 A^2 t^2 + \frac{2}{e t^2}] e^{-\sqrt{2} \uo t +  A t^2}.
\end{equation}
Inserting (\ref{var02}) and (\ref{var3}) into (\ref{var01}) gives the  claim.
\qed

\subsection{Proof of Lemma \ref{lemma:fg}}
For short, we drop $t$ from the functions whenever there is no confusion.   For the first claim,   by direct calculations, we have:
\[
\uo(g,t) - \uo(f,t)  =   \frac{\frac{d}{dt}|f|}{|f|}  -   \frac{\frac{d}{dt}|g|}{|g|} =    I +   II,
\]
where   $I  = (1 - \frac{|g|^2}{|f|^2}) \cdot \uo(g,t)$,
$
 II = \frac{1}{   |f|^2} \cdot    [\mathrm{Re}(g') \cdot  \mathrm{Re}(f  - g)  +  \mathrm{Im}(g') \cdot  \mathrm{Im}(f-g)  +    \mathrm{Re}(f)  \cdot  \mathrm{Re}((f - g)')  +  \mathrm{Im}(f)  \cdot  \mathrm{Im}((f - g)')]$.
Now, firstly,   using triangle inequality,
\[
|I| \leq   \frac{|\uo(g,t)|}{|f|^2} \cdot  \big| |f|^2 - |g|^2 \big|  \leq   \frac{|\uo(g,t)|}{|f|^2} \bigl( (|f| + |g|) |f - g| \bigr);
\]
secondly, using  Cauchy-Schwartz  inequality,
$|\mathrm{Re}(z) \mathrm{Re}(w) + \mathrm{Im}(z) \mathrm{Im}(w)| \leq |z| \cdot |w|$  for any complex numbers $z$ and $w$,   so it follows that
\[
|II| \leq \frac{\sqrt{2}}{ |f|^2}   \cdot [|g'| \cdot |f  - g| + |f| \cdot |(f  - g)'|].
\]
Combining these gives
\[
|\uo(g,t) - \uo(f,t)| \leq \frac{1}{  |f|^2}   \bigl[(    |\uo(g,t)| (|f| + |g|)   +  \sqrt{2} |g'|)   \cdot  |f - g|  + \sqrt{2}  |g| \cdot  |(f - g)'|\bigr].
\]

Consider  the second claim.  By direct calculations,
\[
\so^2(g,t) - \so^2(f,t)  =   I + II,
\]
where  $I  =   (1  - \frac{|g|^2}{|f|^2}) \cdot  \uo(g,t)$,
$II  =  \frac{\sqrt{2}}{t|f|^2} \cdot [(\mathrm{Re}(\omega \bar{g} g') - \mathrm{Re}(\omega \bar{f} f')]$.
Similarly,
\[
|I| \leq  \frac{|\so^2(g,t)|}{|f|^2}  (|f| + |g|)  |f - g|,
\]
\[
|II| \leq  \frac{\sqrt{2}}{t|f|^2} \cdot [ |g'| \cdot |f  -  g|  + |f| \cdot |(f  - g)'|].
\]
Combining these gives
\[
|\so^2(g, t) - \so^2(f,t)|  \leq  \frac{1}{t |f|^2} \cdot \bigl[(|\so^2(g,t)| \cdot t \cdot  (|f| + |g|)  + \sqrt{2} |g'|) \cdot |f - g|  +   \sqrt{2} |f| \cdot |(f - g)'| \bigr].
\]
\qed

\subsection{Proof of Lemma \ref{lemma:stoch}}
Write $t_n = t_n(\gamma)$ for short. Introduce the event
\[
A_n = \{ \max\{ |\varphi_n(t_n) - \varphi(t_n)|,  \;    |\varphi_n'(t_n)  - \varphi'(t_n)|  \}   \leq \log^{3/2}(n)  \}.
\]
Applying Lemma \ref{lemma:var}, $P(A_n^c) \goto 0$, uniformly for all $f \in \Lambda_n(\eps_0, A)$.    To show the claim,  it is sufficient to show that the inequalities hold
over the event $A_n$.  Since the proofs are similar, we only prove the first one.

We claim that (a). $1/\varphi(t_n)| \leq \bar{o}(1)$ over event $A_n$,   (b).
$|\varphi_n(t_n)|  \sim |\varphi(t_n)|$ over event $A_n$, and (c). $|\uo(\varphi; t_n)| \leq \bar{o}(1)$.
Consider (a) and (b).   By $\eps_n \leq \eps_0 < 1/2$ and elementary calculus,  $|r(t)| \leq \eps_n/(1 - \eps_n) \leq \eps_0/(1 - \eps_0)$.
The claim follows from
\[
|\varphi(t)|  \geq (1 - |r(t)|) |\varphi_0(t)| \geq \frac{1 - 2 \eps_0}{1 - \eps_0}  e^{- \uo t_n /\sqrt{2}}.
\]
By the definition of $A_n$,
\[
|\varphi_n(t_n)  - \varphi(t_n)| \leq   \bar{o}(n^{(A \gamma -1)/2}),
\]
and the claims follow.  Consider (c).   By Lemma \ref{lemma:r}, $|\uo(\varphi_n; t)| \leq |\uo| + |r'(t_n)|$.  Write
\[
r'(t) = \frac{\eps_n}{1 - \eps_n}  \int [\omega(u - \uo) + i (\sigma^2 - \so^2)  t]  e^{ \omega (u  - \uo) t + i (\sigma^2  - \so^2) t^2/2} d H_n(u, \sigma).
\]
Since that $\sup_{x > 0} \{ x e^{-x} \}  = 1/e$ and $P_{H_n} (u > \uo,  \sigma^2 \leq A) = 1$,     the claim follows from
\[
|r'(t_n)| \leq   \int [ |u - \uo|  e^{- (u - \uo)t_n/\sqrt{2} }  +  |\sigma^2 - \so^2| t_n]   d H_n(u, \sigma)   \leq  (\sqrt{2}/(et_n)) + A t_n.
\]

Finally, combine (a)-(c) with Lemma \ref{lemma:fg},
\[
|\uo(\varphi_n; t_n)  - \uo(\varphi; t_n)| \leq \bar{o}(1) \cdot [ |\varphi_n(t_n) - \varphi(t_n)| +  |\varphi_n'(t_n) - \varphi'(t_n)|],
\]
and the claim follows. \qed

\subsection{Proof of Lemma \ref{lemma:r}}
For simplicity, drop $t$ from $\varphi(t)$,  $\varphi_0(t)$, and $r(t)$  whenever there is no confusion.
Consider the first claim.  Recalling that  $|\varphi|=|\varphi_0 |\times |1+r|$,
\[
\frac{d}{dt}|\varphi(t)|=(\frac{d}{dt} |\varphi_0|)\cdot  |1+r|+|\varphi_0| \cdot   \frac{d}{dt}|1+r|.
\]
Using the  definition of $\uo(\varphi; t)$ and  Lemma \ref{lemma:varphi0},  it follows from direct calculations that
\begin{equation} \label{lemmar1}
|\uo(\varphi;t) - \uo|  =   \frac{\sqrt{2}}{|1+r(t)|}\frac{d}{dt} |1+r(t)|.
\end{equation}
Moreover,
\begin{equation} \label{lemmar2}
\frac{d}{dt}|1+r(t)|=\frac{r'(t)(1+ \bar{r}(t))+(1+r(t))  \bar{r}'(t)}{2|1+r(t)|}.
\end{equation}
By that $P_{H_n} (u > \uo) = 1$,
\begin{equation} \label{lemmar3}
|r(t)| \leq  \frac{\eps_n}{1 - \eps_n}  \bigl|  \int e^{ - (u - \uo) t /\sqrt{2} + i  ( - (u - \uo)t /\sqrt{2}  + (\sigma^2 - \so^2) t^2/2}  d H_n(u, \sigma) \bigr|   \leq \frac{\eps_n}{1 - \eps_n}.
\end{equation}
Combining(\ref{lemmar1})--(\ref{lemmar3})  gives the claim.

Consider the second claim.
Write
\[
\varphi' =\varphi_0' (1+r)+\varphi_0  r'.
\]
We have  $\bar{\varphi}\varphi' =|1+r|^2\bar{\varphi}_0 \varphi_0' +|\varphi_0|^2(1+\bar{r})r'$, and so
\[
\mathrm{Re}(\omega \bar{\varphi} \varphi)   = |1 + r|^2 \mathrm{Re} (\omega \bar{\varphi}_0 \varphi_0)  + |\varphi_0|^2 \mathrm{Re}(\omega  (1 + \bar{r}) r').
\]
Therefore,
\[
|\sigma^2_0(\phi;t) - \sigma_0^2| \leq \frac{| \mathrm{Re} (\omega  (1+\bar{r}(t)) r' )| }{t|1+r(t)|^2} \leq  C(\eps_0)    |r'(t)|/t,
\]
and the claim follows directly.    \qed

\subsection{Proof of Theorem \ref{thm:eps1}   }
Write for short  $t_n=t_n(\gamma)$ and $\eps_n(\cdot; t_n) = \eps_n(\cdot; t_n, \uo, \so^2)$.  By triangle inequality,
\[
|\epsilon_n(\varphi_n;t_n) -   \epsilon_n| \leq |\epsilon_n(\varphi_n;t_n)-   \epsilon_n(\varphi;t_n)|+|\epsilon_n(\varphi;t_n)-\epsilon_n|.
\]
Compare this with the desired claim. It is sufficient to show that  $E[|\epsilon_n(\varphi_n;t_n)-   \epsilon_n(\varphi;t_n)|^2]   \leq    n^{A \gamma -1}$ and $|\epsilon_n(\varphi;t_n)/\eps_n - 1|$ tends to $0$ in a speed that does not depends on $\eps_n$ and $A_n$.

Consider the first term first.  By the definition of the functional $\eps_n(\cdot; t_n)$ (i.e. (\ref{Defineeps})),
\[
|\epsilon_n(\varphi_n;t_n)  -   \epsilon_n(\varphi;t_n)|  \leq  |e^{-\omega u_0 t_n-i\sigma_0^2t_n^2/2}(\varphi_n(t_n)-\varphi(t_n))|
 \leq e^{\frac{u_0 t_n}{\sqrt{2}}}|\varphi_n(t_n)-\varphi(t_n)|.
\]
At the same time, by the definitions of $\varphi_n(\cdot)$ and $\varphi(\cdot)$ and elementary calculus,
\[
E| \varphi_n(t_n) - \varphi(t_n)|^2    =    \frac{1}{n}  \mathrm{Var} ( e^{ \omega t_n X_1})  \leq \frac{1}{n} E[e^{-\sqrt{2} t_n X_1}].
\]
Combining these gives,
\[
E(|\epsilon_n(\varphi_n;t_n)  -   \epsilon_n(\varphi;t_n)|^2)  \leq   \frac{1}{n} e^{\sqrt{2} t_n \uo}  E[e^{-\sqrt{2} t_n X_1}],
\]
where by direct calculations and the assumptions of $P_{H_n}(u>u_0, \sigma^2 \leq A)=1$ and $\so^2 \leq A$,     the last term is no greater than
\[
\frac{1}{n}\left\{(1-\epsilon_n) e^{t_n^2 \sigma_0^2}+\epsilon_n \int e^{-\sqrt{2} t_n(u-u_0)+t_n^2\sigma^2}d H_n(u,\sigma)\right\} \leq n^{A \gamma -1}.
\]
Combining these gives the first claim.

Consider the second claim.  Recall that  $\epsilon_n=\epsilon_n(\varphi_0;t_n)$,  that $\varphi(t)  = \varphi_0(t)( 1 + r(t))$ (see (\ref{Definer})), and that $\varphi_0(t) = (1 - \eps_n) e^{-\omega u_0 t_n-i\sigma_0^2t_n^2/2}$.
By the definition of the functional $\eps_n(\cdot, t_n)$,
\[
\eps_n(\varphi; t_n)   - \eps_n =  e^{-\omega u_0 t_n-i\sigma_0^2t_n^2/2}(\varphi(t_n)-\varphi_0(t_n))  = (1 - \eps_n)  r(t_n).
\]
It then follows from   the definition of $r(\cdot)$  that
\begin{equation} \label{rterm0}
|\eps_n(\varphi; t_n)   - \eps_n|   \leq   |(1 - \eps_n) r(t_n)|  =  \epsilon_n\left|\int e^{\omega (u-u_0)t_n+i(\sigma^2-\sigma_0^2)t_n^2/2} d H_n(u,\sigma)\right|.
\end{equation}

Suppose condition (b) holds.  Then $P_{H_n}(u>u_0+\delta_n)=1$,  where $\delta_n$ satisfies (\ref{Definedelta}) with some constant $c_0 > 0$.      It follows from (\ref{rterm0}) and elementary calculus that  as $n \goto \infty$,
\begin{equation} \label{rterm}
| \eps_n(\varphi; t_n)/\eps_n   - 1 |
\leq    \int e^{-\frac{(u-u_0)t}{\sqrt{2}}} d H_n(u,\sigma)     \goto 0.
\end{equation}

Suppose condition (c) holds.
Let $\delta =u-u_0, \kappa =\sigma^2-\sigma_0^2$.  By the definitions of $g(\kappa | \delta)$ and $g(\delta)$ and  elementary Fourier analysis,
\[
\int e^{\omega (u-u_0)t_n+i(\sigma^2-\sigma_0^2)t_n^2/2} d H_n(u,\sigma)  =   \int e^{\omega \delta t_n+i\kappa t_n^2/2}  g(\kappa|\delta) h(\delta)d\kappa d\delta
 =  \int e^{\omega \delta t_n} g^{FT}(t_n^2/2;\delta) h(\delta) d\delta.
\]
By the assumptions,   $P_{H_n}{(\delta>0)}=1$ and  $g^{FT}(t) \leq C(1 + |t|)^{-\alpha}$, so
\begin{eqnarray*}
\big| \int e^{\omega \delta t_n} g^{FT}(t_n^2/2;\delta) h(\delta) d\delta  \bigr|  \leq  \int e^{-\frac{\delta t_n}{\sqrt{2}}} |g^{FT}(t_n^2/2;\delta)| h(\delta) d\delta
\leq C(1+t_n^2/2)^{-k} \longrightarrow 0,
\end{eqnarray*}
Combining these with (\ref{rterm0}) gives the claim.  \qed

\subsection{Proof of Theorem \ref{thm:eps2}}
Write for short  $t_n=t_n(\gamma)$, $\hat{u}_0 =u_0(\varphi_n;t_n)$, and $\hat{\sigma}_0^2 =\sigma_0^2(\varphi_n;t_n)$.  By the definitions of $\eps_n(\cdot;  t, u, \sigma)$,
\[
|\eps_n(\varphi_n; t_n, \hat{u}_0, \hat{\sigma}_0^2)  -  \eps_n(\varphi_n; t_n, \uo, \so^2)|     \leq | e^{ - \omega (\hat{u}_0 - \uo) t_n - i (\hat{\sigma}_0^2 - \so^2) t_n^2} - 1| \cdot | e^{- \omega \uo t_n  -  i \so^2 t_n^2} \varphi_n(t_n)|,
\]
where we note that by the definition of the functional  $\eps_n(\cdot; t, u, \sigma)$,   the last term $\leq  1 + |\eps_n(\varphi_n; t_n, \uo, \so^2)|$.
By Lemmas \ref{lemma:stoch}--\ref{lemma:r}, except for a small probability that tends to $0$ as $n$ tends to $\infty$,
\[
|\hat{u}_0 - \uo|  t_n  \leq   t_n |r'(t_n)| + \bar{o}(n^{(A \gamma -1)/2}), \qquad |\hat{\sigma}_0^2 - \so^2| t_n^2  \leq  t_n |r'(t_n)| + \bar{o}(n^{A \gamma  - 1)/2}).
\]
At the same time,  by Theorem \ref{thm:eps1}, except for  a small probability that tends to $0$ as $n$ tends to $\infty$,
\[
\eps_n(\varphi_n; t_n, \uo, \so^2)  \sim \eps_n.
\]
Combine these, as $n$ tends to $\infty$,   except for  a small probability that tends to $0$.
\[
|\eps_n(\varphi_n; t_n, \hat{u}_0, \hat{\sigma}_0^2)  -  \eps_n(\varphi_n; t_n, \uo, \so^2)|   \leq t_n |r'(t_n)|,
\]
which, by Lemmas \ref{lemma:exp1}--\ref{lemma:exp2}, tends to $0$. This concludes the proof.    \qed


\end{document}